\documentstyle[12pt]{article}
\newtheorem{theorem}{Theorem}[section]
\newtheorem{proposition}[theorem]{Proposition}

\newtheorem{lemma}[theorem]{Lemma}
\newtheorem{corollary}[theorem]{Corollary}

\newtheorem{definition}[theorem]{Definition}

\font\sym=msbm10 scaled \magstep1
\newcommand{\IZ}{\hbox{\sym \char '132}}
\newcommand{\IP}{\hbox{\sym \char '120}}

\newcommand{\ka}{{\cal A}}
\newcommand{\kb}{{\cal B}}
\newcommand{\kd}{{\cal D}}
\newcommand{\ke}{{\cal E}}
\newcommand{\kl}{{\cal L}}
\newcommand{\km}{{\cal M}}
\newcommand{\ko}{{\cal O}}
\newcommand{\kt}{{\cal T}}
\newcommand{\kz}{{\cal Z}}
\newcommand{\kx}{{\cal X}}

\newcommand{\Pic}{{\rm Pic }}
\newcommand{\Hom}{{\rm Hom }}
\newcommand{\Ext}{{\rm Ext }}

\newcommand{\Ker}{{\rm Ker }}
\newcommand{\Coker}{{\rm Coker }}
\newcommand{\qed}{\hspace*{\fill}\hbox{$\Box$}}

\newcommand{\lra}{\longrightarrow}
\newcommand{\sesq}[5]{0\lra#1\stackrel{#2}{\lra}#3\stackrel{#4}{\lra}#5\lra0}

\begin{document}
\title{Hodge numbers of moduli spaces of stable bundles on K3 surfaces}
\author{ L.\ G\"ottsche~~~~~ D.\ Huybrechts}
\date{}
\maketitle

Very few higher dimensional symplectic manifolds are known. Beauville
has shown that the Hilbert schemes of points on a K3 surface are
symplectic \cite{B1}.
Other examples are provided by the moduli spaces of stable vector bundles on
these surfaces \cite{Mu1}. Very often both spaces are closely related
by a birational correspondence, but in general
they are not isomorphic (\cite{Mu2}, p. 167). Beauville also proved that all
except one of the
deformation directions of the Hilbert scheme are obtained by deforming the
underlying K3 surface. It is natural to ask if in fact the moduli spaces
are obtained by deforming the Hilbert scheme in this extra direction.
The aim of this paper is to show that
for some of the moduli spaces an important class of deformation invariants,
namely the Hodge numbers, coincide with those of an appropriate Hilbert scheme.
The following theorem is proven:

{\bf Theorem:} {\it Let $X$ be a K3 surface, $L$ a primitive big and
nef line bundle and $H$ a generic polarization.
If $\overline{M}_H(L,c_2)$ denotes the moduli space
of rank two semi-stable torsion-free sheaves and
$\dim\overline{M}_H(L,c_2)>8$ then its Hodge numbers
coincide with the Hodge numbers
of the Hilbert scheme of $l:=2c_2-\frac{L^2}{2}-3$ points,}
i.e.$$h^{p,q}(\overline{M}_H(L,c_2))=h^{p,q}(Hilb^l(X)).$$

By a generic polarization we mean a polarization which does not lie on any
wall, i.e. any $H-$semi-stable sheaf is stable.
Note that by the smoothness criterion \cite{Mu1} the moduli space
$\overline{M}_H(L,c_2)$ is smooth and projective of dimension $4c_2-L^2-6$
if it is non-empty.
Moreover, by dimension counting (cf. \cite{T}, lemma 3.1, \cite{O'G}, prop.
7.2.1)
the moduli space of $\mu-$stable vector
bundles is dense in $\overline{M}_H(L,c_2)$.
We would like to mention that the Hodge numbers of the Hilbert scheme of
points on a surface can be expressed in terms of the Hodge numbers of
the surface \cite{G1,Ch}. In fact, even the Hodge structure is known.
We would like to compare the Hodge structure of the moduli spaces
and of the Hilbert schemes, but even in the special case dealt with in section
1 we cannot prove that they coincide.

Our work was motivated by a talk of J. Le Potier in Lambrecht in May 1994.
In this talk he explained how to use moduli spaces of coherent systems
(or framed modules) to compute the $SU(2)-$Donaldson polynomials for
K3 surfaces. This was done before by O'Grady by other methods \cite{O'G1}.
We tried to use Le Potier's approach to compute the $SO(3)-$polynomials
which turned out to be even easier. (In the meantime the computation
of them using
O'Grady's method
has appeared in \cite{Na}. Therefore our computation is not
included.) Then we realized that the approach also
works in the context of this paper.

{\bf Notations:} $X$ will always be a projective K3 surface.\\
$P_E$ - the Hilbert polynomial of a sheaf $E$ (we suppress the fixed
polarization $H$ in the notation).\\
$M_H(L,c_2)$ - the moduli space of rank two vector bundles
with determinant $L\in \Pic(X)$ and second Chern class $c_2$ which are
$\mu-$stable with respect to a polarization $H$.\\
$\overline M_H(L,c_2)$ - the moduli space of rank two torsion-free sheaves
with determinant $L$ and second Chern class $c_2$ which are semi-stable with
respect to $H$.\\
$M_H(L,c_2,D,\delta)$ - moduli space of framed modules $(E,\alpha:E\to D)$,
$E$ locally free with the given invariants, which are $\mu-$stable with respect
to $\delta$ and $H$.\\
$\overline M_H(L,c_2,D,\delta)$ - the moduli space of all semi-stable framed
modules.\\
If $H=L$ we drop $H$ in
the notation.
\section{A special case}
In this section we prove the theorem in the case that $\Pic(X)=\IZ\cdot
L$ and $c_2=\frac{L^2}{2}+3$.
{\small\subsection{ The birational correspondence to the Hilbert scheme}}

Throughout this section we will assume that the Picard group is generated by
an ample line bundle $L$, i.e. $\Pic(X)=\IZ \cdot L$.
Under this assumption a
torsion-free sheaf with determinant $L$ is $\mu-$stable if and only if it is
$\mu-$semi-stable.

For the convenience of the reader we recall the stability condition for framed
modules (\cite{HL}):
Let $\delta=\delta_1\cdot n+\delta_0$, $D\in \Pic(X)$ and $E$ be a torsion-free
rank two sheaf. A framed module $(E,\alpha)$ consists of $E$ and a
non-trivial homomorphism
$\alpha:E\to D$. It is (semi-)stable if
$P_{\Ker(\alpha)}(\leq)P_E/2-\delta/2$ and for all rank one subsheaves
$M\subset E$ the inequality $P_M(\leq)P_E/2+\delta/2$ holds.\\
In fact, any semi-stable framed module is torsion-free. In \cite{HL} it was
shown
that there exists a coarse projective moduli space of semi-stable
framed modules.

\begin{lemma}
Let $D=L$ and $0<\delta_1<L^2$. Then a framed module $(E,\alpha)$ is
$\mu-$stable if and only if $E$ is $\mu-$stable.  The moduli space
$\overline M_H(L,c_2,D,\delta)$ is independent of the specific $\delta$ in
this range.
\end{lemma}
{\it Proof:} Let $(E,\alpha)$ be semi-stable, then $\mu(M)\leq
L^2/2+\delta_1/2$
for all $M=L^{\otimes n}\otimes I_Z\subset E$. Thus $nL^2\leq
L^2/2+\delta_1<L^2$, i.e.
$n<1$. Hence $E$ is $\mu-$stable. If $E$ is $\mu-$stable, then
$(E,\alpha)$ is $\mu-$stable if $\mu(\Ker(\alpha))<L^2/2-\delta_1/2$. But
writing $\Ker(\alpha)=L^{\otimes n}\otimes I_Z$ and using the stabilty of $E$
we conclude
$n<1$. Hence $\mu(\Ker(\alpha))\leq0<L^2/2-\delta_1/2$. The second statement
follows immediately.\qed

Henceforth $\delta$ is chosen as in the lemma. Note that as for sheaves
also for framed modules $\mu-$stability is equivalent to $\mu-$semi-stability.
It can also be shown that both moduli spaces $ M(L,c_2,L,\delta)$ and
$\overline M(L,c_2)$ are fine.
The universality property of the moduli space induces a morphism
$$\overline\varphi:\overline M(L,c_2,L,\delta)\lra\overline M(L,c_2).$$
%

Note that by the previous lemma the fibre of
$\overline\varphi$ over $[E]$ is isomorphic
to $\IP(\Hom(E,L))$.
\begin{lemma}\label{phisur} If $c_2\leq L^2/2+3$, then $\varphi$ is surjective.
\end{lemma}
{\it Proof:} It suffices to show that for a $\mu-$stable vector bundle
$E$ there is always a non-trivial homomorphism $E\to L$.
Since $\Hom(E,L)=H^0(X,E)$
and $H^2(X,E)\cong H^0(X,E^*)=0$ by the stability of $E$,
the Riemann-Roch-Hirzebruch formula $\chi(E)=L^2/2-c_2+4$ shows that under
the assumption $h^0(X,E)>0$.\qed

\begin{lemma}\label{defN} Let $N(L,c_2,L,\delta)$ be the set of all
$(E,\alpha)\in\overline {M}(L,c_2,L,\delta)$ such that
$\Ker(\alpha)\cong\ko_X$. It is a closed subset, which contains all stable
pairs
$(E,\alpha)$ with $E$ locally free.
\end{lemma}
{\it Proof:} If $E$ is locally free, then $\Ker(\alpha)$ has to be locally
free. By stability it is thus isomorphic to $\ko_X$.
 $N(L,c_2,L,\delta)$ is closed, since  $\Ker(\alpha)\cong\ko_X$ if and only if
$length(\Coker(\alpha))=c_2$, i.e. if the length is maximal;
this is a closed condition.\qed

Under the condition of \ref{phisur} we have a surjective
morphism$$\varphi:{N(L,c_2,L,\delta)}\lra\overline M(L,c_2).$$
By \ref{defN} any framed module $(E,\alpha)\in {N(L,c_2,L,\delta)}$
sits in
an extension
$$\sesq{\ko_X}{}{E}{\alpha}{I_Z\otimes L},$$
where $I_Z$ is the ideal sheaf of a codimension two cycle of length $c_2$.
Thus we can define a morphism
$$\psi: {N(L,c_2,L,\delta)}\lra Hilb^{c_2}(X)$$
by mapping $(E,\alpha)$ to $[\Coker(\alpha)]$.

\begin{lemma}, If $c_2\geq L^2/2+3$, $\psi$ is surjective.
\end{lemma}
{\it Proof:} It is enough to show that $\Ext^1(I_Z\otimes L,\ko_X)\not=0$
for all $Z\in Hilb^{c_2}(X)$. By the assumption $h^0(X,L|_Z)-h^0(X,L)\geq1$.
Thus $h^1(L\otimes I_Z)\geq 1$. Now use
$\Ext^1(I_Z\otimes L,\ko_X)\cong H^1(X,I_Z\otimes L)^*$.\qed

We have seen that any $(E,\alpha)\in {N(L,c,L,\delta)}$ induces
an exact sequence
$$\sesq{\ko_X}{}{E}{\alpha}{I_Z\otimes L}.$$
Conversely, any section $s\in H^0(X,E)$ of $E\in\overline M(L,c)$
gives a homomorphism $\alpha:E\to L$ with $\Ker(\alpha)\cong\ko_X$.
Thus the fibre of $\varphi: {N(L,c,L,\delta)}
\to \overline M(L,c)$ over $[E]$ is isomorphic to $\IP(H^0(X,E))$. In fact,
$ {N(L,c,L,\delta)}$ can be identified with Le Potier's moduli space
of coherent systems of rank one \cite{LP}.

The picture we get in the case $c_2=c:=L^2/2+3$ is described by the following
diagram.
$$\begin{array}{rlcrl}
&& {N(L,c,L,\delta)}&&\\
&\varphi\swarrow&&\searrow\psi&\\
{}~~~~\overline{M}(L,c)&&&&Hilb^{c}(X)~~~~\\
\end{array}$$
Both morphisms $\varphi$ and $\psi$ are birational. This is
due to the fact that for the generic $[Z]\in Hilb^{c}(X)$ the
restriction map $H^0(X,L)\to H^0(X,L_Z)$ is injective and hence $h^1(X,I_Z
\otimes L)=1$. This shows that $\psi$ is generically an isomorphism.
Since the fibres of $\varphi$ are connected and both spaces are of the same
dimension, also $\varphi$ is birational. Note that in particular the
moduli space $\overline{M}(L,c)$ is irreducible.\\
Results about birationality of certain moduli spaces and corresponding
Hilbert schemes have been known for some time, e.g. Zuo has shown that
$M_H(\ko_X,n^2H^2+3)$ is birational to $Hilb^{2n^2H^2+3}(X)$ (cf. \cite{Z}).
The moduli spaces of framed modules make this  relation more explicit. They
are used in the next section to show that the Hodge numbers of the moduli space
and the Hilbert scheme coincide.\\

{\small\subsection{Comparison of the Hodge numbers}}

First, we recall the notion of virtual Hodge polynomials \cite{D}, \cite{Ch}.\\
For any quasi-projective variety $X$ there exists a polynomial $e(X,x,y)$
with the following properties:\\
{\it i)} If $X$ is smooth and projective then
$$e(X,x,y)=h(X,-x,-y):=\sum_{p,q}(-1)^{p+q}h^{p,q}(X)x^py^q.$$
{\it ii)} If $Y\subset X$ is Zariski closed and $U$ its complement then
$$e(X,x,y)=e(Y,x,y)+e(U,x,y).$$
{\it iii)} If $X\to Y$ is a Zariski locally trivial fibre bundle with fibre $F$
then $$e(X,x,y)=e(Y,x,y)\cdot e(F,x,y).$$
{\it iv)} If $X\to Y$ is a bijective morphism then $e(X,x,y)=e(Y,x,y)$.

In particular, if
$$\begin{array}{ccccc}
&&Z&&\\
&\swarrow&&\searrow\\
X&&&&Y\\
\end{array}$$
is a diagram of quasi-projective varieties,
where $Z\to X$ and $Z\to Y$ admit a bijective morphism to a $\IP_n-$ bundle
over $X$, resp. $Y$, then
$e(X,x,y)\cdot e(\IP_n,x,y)=e(Z,x,y)=e(Y,x,y)\cdot e(\IP_n,x,y)$.
Hence $e(X,x,y)=e(Y,x,y)$.\\
The idea to prove that $\overline{M}(L,c)$ and $Hilb^{c}(X)$, with
$ c:=\frac{L^2}{2}+3$, have the same Hilbert polynomial is to stratify both by
locally
closed subsets $\overline{M}(L,c)_k$ and $Hilb^c(X)_k$
such that the birational correspondence given by the moduli space of
framed modules induces $\IP_{k-1}-$bundles ${N(L,c,L,\delta)}_k \to
\overline{M}(L,c)_k$ and ${N(L,c,L,\delta)}_k \to Hilb^c(X)_k$.
One concludes $e(\overline{M}(L,c)_k,x,y)=e(Hilb^c(X)_k,x,y)$ and hence
$e(\overline{M}(L,c),x,y)=\sum_k e(\overline{M}(L,c)_k,x,y)=\sum_k
e(Hilb^c(X)_k,x,y)=e(Hilb^c(X),x,y)$.

We first define the stratification.
\begin{definition}
$Hilb^c(X)_k:=\{[Z]\in Hilb^c(X)~|~h^1(X,I_Z\otimes L)=k\}$\\
${N(L,c,L,\delta)}_k:=\psi^{-1}(Hilb^c(X)_k)$\\
$\overline{M}(L,c)_k:=\varphi({N(L,c,L,\delta)}_k)$
\end{definition}
Using the universal subscheme $\kz\subset X\times Hilb^c(X)$ with the
two projections $p$ and $q$ to $X$ and $Hilb^c(X)$, resp., and the
semi-continuity applied to the sheaf $I_\kz\otimes p^*(L)$ and the projection
$q$ it is easy to see that this defines a stratification into locally
closed subschemes. All strata are given the reduced induced structure.

We want to show that both morphisms
$${N(L,c,L,\delta)}_k\to Hilb^c(X)_k$$ and
$${N(L,c,L,\delta)}_k\to \overline{M}(L,c)_k$$ admit a bijective morphism to
a
$\IP_{k-1}-$bundle over the base. In fact, they are $\IP_{k-1}-$bundles, but by
property
{\it iv)} of the virtual Hodge polynomials we only need the bijectivity.
\begin{definition} Let $\ka_k:=\pi_*(\ke_k)$,
where
$\pi:{N(L,c,L,\delta)}_k\times X\to {N(L,c,L,\delta)}_k$
denotes the projection and $\ke_k$ is the restriction of the universal sheaf
$\ke$, and
 let $\kb_k:={\cal E}xt^1_q((I_\kz)_k\otimes p^*(L),\ko_X)$ be
the relative Ext-sheaf, where $(I_\kz)_k$ denotes the restriction of
$I_\kz$ to $ Hilb^c(X)_k\times X$.
\end{definition}

\begin{lemma} $\ka_k$ and $\kb_k$ are locally free
sheaves on $\overline{M}(L,c)_k$ and $Hilb^c(X)_k$, resp.,
and compatible with base change,
i.e. $\ka_k([E])\cong H^0(X,E)$ and $\kb_k([Z])\cong \Ext^1(I_Z\otimes
L,\ko_X)$.
\end{lemma}
{\it Proof:} By definition and using Serre-duality
we see that $Hilb^c(X)_k=\{Z|\dim\Ext^1(I_Z\otimes L,\ko_X)=k\}$ and that
it is reduced. Thus the claim for $\kb_k$ follows immediately from
the  base change theorem for global Ext-groups \cite{BPS}. In order to
prove the assertion for $\ka_k$ it suffices to
show that $\overline M(L,c)_k=\{E\,|\,h^0(X,E)=k\}$.
Consider the exact sequence
$$\sesq{\ko_X}{}{E}{}{I_Z\otimes L}.$$
Then $E\in \overline M(L,c)_k$ if and only if $h^1(X,I_Z\otimes L)=k$ if
and only if $h^0(X,I_Z\otimes L)=k-1$ if and only if $h^0(X,E)=k$.\qed

The kernel of the universal framed module on ${N(L,c,L,\delta)}\times X$
restricted
to ${N(L,c,L,\delta)}_k$
induces a morphism to $\IP(\ka_k)$ which is obviously
bijective.
Analogously, by the universality of $\IP(\kb_k)$ (cf. \cite{La}) the
universal framed module over ${N(L,c,L,\delta)}\times X$
completed to an exact sequence and restricted to the stratum
induces a bijective morphism
of ${N(L,c,L,\delta)}_k$ to $\IP(\kb_k)$.

We summarize:
\begin{proposition}
If $X$ is a K3 surface with $\Pic(X)=\IZ\cdot L$, $L$ ample
and $c_2=L^2/2+3$, then $h^{p,q}(\overline{M}(L,c_2))=h^{p,q}(Hilb^{c_2}(X))$.
\qed
\end{proposition}

Both manifolds $\overline{M}(L,c_2)$ and $Hilb^{c_2}(X)$
are symplectic. One might conjecture that in general two birational
symplectic manifolds have the same Hodge numbers or even isomorphic Hodge
structures,
but we don't know how to prove this.


\section{The general case}

By deforming the underlying K3 surface the proof of the theorem is
reduced to the case considered in section 1.

\subsection{Deformation of K3 surfaces}

The following statements about the existence of certain
deformations of a given K3 surface will be needed.

{\bf 2.1.1} {\it  Let $X$ be a K3 surface, $L\in\Pic(X)$ a primitive nef and
big
line bundle. Then there exists a smooth connected family $\kx\lra S$
of K3 surfaces and a line bundle $\kl$ on $\kx$ such that:

$\cdot$ $\kx_0\cong X$ and $\kl_0\cong L$.

$\cdot$ $\Pic(\kx_t)=\IZ\cdot\kl_t$ for all $t\not=0$ ($\kl_t$ is automatically
ample).}

{\it Proof:} The moduli space of primitive pseudo-polarized K3 surfaces is
irreducible (\cite{B2}). Since any even lattice of index $(1,\rho-1)$
with $\rho\leq10$ can be realized as a Picard group of a K3 surface
(\cite{Ni},\cite{Mor}) the generic pseudo-polarized K3 surface has Picard group
$\IZ$.\qed

{\bf 2.1.2} {\it  Let $X$ be a K3 surface whose Picard group is generated by an
ample
line bundle $L$, i.e. $\Pic(X)=\IZ\cdot L$.
Furthermore, let $d\geq 5$ be an integer.
Then there exists a smooth connected family $\kx\lra S$ of K3 surfaces
and a line bundle $\kl$ on $\kx$
such that:

$\cdot$ $(\kx_{t_0},\kl_{t_0})\cong(X,L)$ for some point $t_0\in
S\setminus\{0\}$.

$\cdot$ $\Pic(\kx_t)\cong \IZ\cdot\kl_t$ for all $t\not=0$.

$\cdot$ $\Pic(\kx_0)=\IZ\cdot\kl_0\oplus\IZ \cdot D$, where $D$ is represented
by a smooth
rational curve, both line bundles $\kl_0$ and $\kl_0(2D)$ are ample
and primitive
and the intersection matrix is}
$$\left(\begin{array}{cc}L^2&d\\
d&-2\
\end{array}\right)$$

{\it Proof:} Again we use the irreducibility of the moduli space of
primitive polarized K3 surfaces. The existence of a triple
$(\kx_0,\kl_0,D)$ with ample $\kl_0$, smooth rational $D$ and the
given intersection form was shown by Oguiso \cite{Og}. It remains to show that
$\kl_0(2D)$ is ample. Obviously, $\kl_0(2D)$ is big
and for any irreducible curve $C\not=D$ the strict inequality $(\kl_0(2D)).C>0$
holds. The assumption on $d$ implies $(\kl_0(2D)).D>0$.
Note that the extra assumption $L^2\geq4$ in \cite{Og} is only needed
for the very ampleness of $\kl_0$ which we will not use.\qed

{\bf 2.1.3(a)} {\it  Let $X$ be a K3 surface whose Picard group is generated by
an ample
line bundle $L$, i.e. $\Pic(X)=\IZ\cdot L$. If $L^2>2$
there exists a smooth connected family $\kx\lra S$ of K3 surfaces
and a line bundle $\kl$ on $\kx$
such that:

$\cdot$ $(\kx_{t_0},\kl_{t_0})\cong(X,L)$ for some point $t_0\in
S\setminus\{0\}$.

$\cdot$ $\Pic(\kx_t)\cong \IZ\cdot\kl_t$ for all $t\not=0$.

$\cdot$ $\Pic(\kx_0)=\IZ\cdot\kl_0\oplus\IZ \cdot D$, where both line bundles
$\kl_0$ and $\kl_0(2D)$ are ample
and primitive
and the intersection matrix is}
$$\left(\begin{array}{cc}L^2&1\\
1&0\
\end{array}\right)$$

{\bf 2.1.3(b)} {\it If we assume that $L^2>6$ we have the same result as in (a)
with ``$\kl_0(2D)$ is ample'' replaced by ``$\kl_0(-2D)$ is ample''.}

{\it Proof:} For both parts we need to prove the existence of a triple
$(X_0,H,D)$ with ample and primitive $H$ and $H(2D)$,
such that $D^2=0$, $H.D=1$ and $H^2=2n>2$ for given $n$. By the results of
Nikulin
we can find a K3 surface with this intersection form. It remains to show that
$H$ and $H(2D)$ can be chosen ample.
We can assume that $H\in{\cal C}^+$, i.e. $H$ is in the positive
component of the positive cone (if necessary change $(H,D)$ to
$(-H,-D)$). We check that $H$ is not orthogonal
to any (-2) class, i.e. for any $\delta:=aH+bD$ ($a,b\in\IZ$) with
$\delta^2=a^2H^2+2ab=-2$ we have $H.\delta\not=0$.
If $H$ were orthogonal to $\delta$ this would
imply that $aH^2+b=0$. Hence $-a^2H^2=-2$ which contradicts $H^2>2$.
Thus $H$ is contained in a chamber. Since the Weyl group
$W_{X_0}$, which is generated
by the reflection on the walls, acts transitively on the set of chambers, we
find $\sigma\in W_{X_0}$ such that $\sigma(H)$ is contained in the chamber
$\{w\in{\cal C}^+|w\delta>0 {\rm ~for~all~effective~(-2)~classes~}\delta\}$.
Applying $\sigma$ to $(H,D)$ we can in fact assume that $H$ is
contained in this chamber.
On a K3 surface the effective divisors are generated by the
effective (-2) classes and $\overline{{\cal C}^+}\setminus\{0\}$. On both sets
$H$ is positive. Thus $H$ is ample.
In order to prove that also $H(2D)$ is ample we show that $D$ is effective
and irreducible.
This follows from the Riemann-Roch-Hirzebruch formula $\chi(\ko(D))=2$,
which implies $D$ or $-D$ effective, and $H.D=1$. Thus $C.D\geq0$
for any curve $C$. Thus $H(2D).C>0$.
Since $H(2D)$ is big we conclude that $H(2D)$ is ample.
To prove (a) we choose $H^2:=L^2$ and use the irreducibility
of the moduli space to show that $(X,L)$ degenerates to $(X_0,H)$.
Defining $\kl_0:=H$ this proves (a).
In order to prove (b) we fix $(H(2D))^2:=L^2$ and let $(X,L)$ degenerate
to $(X_0,H(2D))$. The assumption on $H$ translates to $L^2>6$.
With $\kl_0:=H(2D)$ we obtain (b).
\qed

\subsection{Deformation of the moduli space}

We start out with the following
\begin{lemma}\label{def} Let $E$ be a simple vector bundle on a K3 surface such
that
$L:=det(E)$ is big. The joint deformations of $E$ and $X$ are unobstructed,
i.e. $Def(E,X)$ is smooth. Moreover, $Def(E,X)\to Def(X)$ and
$Def(L,X)\to Def(X)$ have the same image.
\end{lemma}
{\it Proof:} The infinitesimal deformations
of a bundle $E$ together with its underlying manifold $X$ are paramatrized
by $H^1(X,\kd_0^1(E))$, where $\kd_0^1(E)$ is the sheaf of differential
operators of order $\leq1$ with scalar symbol. The obstructions are elements
in the second cohomology of this sheaf. Using the symbol map we have a short
exact sequence
$$\sesq{\ke nd(E)}{}{\kd_0^1(E)}{}{\kt_X}.$$
Its long exact cohomology sequence
$$H^1(X,\kd_0^1(E))\to H^1(X,\kt_X)\to H^2(X,\ke nd(E))\to H^2(X,\kd_0^1(E))\to
0$$
compares the deformations of $E$, $X$, and $(E,X)$.
In particular, if $E$ is simple the trace homomorphism $H^2(X,\ke nd(E))\to
H^2(X,\ko_X)$ is bijective and the composition with the boundary map
$H^1(X,\kt_X)\to H^2(X,\ke nd(E))$ is the cup-product with $c_1(E)$.
Since there is exactly one direction in which a big and nef line bundle
$L$ cannot be deformed with $X$ the cup-product with $c_1(L)=c_1(E)$ is
surjective.
Thus $H^1(X,\kd_0^1(E))\to H^1(X,\kt_X)$ is onto the algebraic
deformations of $X$ and $H^2(X,\kd_0^1(E))$
vanishes.\qed

The following lemma will be needed for the next proposition.
Its proof is quite similar to what we will use to prove the theorem.
\begin{lemma}\label{irr} If $\Pic(X)=\IZ\cdot L$, then $\overline M(L,c_2)$ is
irreducible for $\dim\overline M(L,c_2)=4c_2-L^2-6>8$.
\end{lemma}
{\it Proof:}
{\it 1st step:} First, we show that $\overline M_H(L,\frac{L^2}{2}+3)$
is irreducible whenever $L$ is an ample line bundle on a K3 surface.\\
By a result of \cite{Q} the moduli spaces $\overline M_H(L,\frac{L^2}{2}+3)$
and
$\overline M_L(L,\frac{L^2}{2}+3)$ are birational.
In particular, the number of irreducible components is the same.
We consider a deformation as in 2.1.1.
The corresponding family of moduli spaces
$\overline M_{\kl_t}(\kl_t,\frac{L^2}{2}+3)$ is proper and by lemma
\ref{def} every stable bundle on $X$ can be deformed to a stable bundle
on any nearby fibre. This shows that $\overline
M_{\kl_0}(\kl_0,\frac{L^2}{2}+3)$
has as many irreducible components as $\overline
M_{\kl_{t\not=0}}(\kl_{t\not=0},\frac{L^2}{2}+3)$, which is irreducible.\\
{\it 2nd step:} Assume $e:=c_2-\frac{L^2}{2}-3>0$ and $L^2>2$. We apply
2.1.3(a).
By the same arguments as above we obtain that the number of irreducible
components of $\overline M(L,c_2)$ is at most the number of irreducible
components of $\overline M_{\kl_0}(\kl_0,c_2)$. Again using \cite{Q}
we know that $\overline M_{\kl_0}(\kl_0,c_2)$ is birational to
$\overline M_{\kl_0(2D)}(\kl_0,c_2)$. The $\mu$-stable part of the latter
is isomorphic to the $\mu$-stable part of
$\overline M_{\kl_0(2D)}(\kl_0(2D),c_2+1)$. We have $(\kl_0(2D))^2=L^2+4>2$
and $c_2+1-\frac{(\kl_0(2D))^2}{2}-3=e-1$. Therefore
we obtain by induction over $e$
and step 1 that $\overline M_{\kl_0(2D)}(\kl_0(2D),c_2+1)$ is irreducible.
Since the locally free $\mu$-stable
sheaves are dense in the moduli spaces, this accomplishes the proof in
this case.\\
{\it 3rd step:} Here we assume that $e:=c_2-\frac{L^2}{2}-3<0$. By assumption
$4c_2-L^2-6\geq10$. Hence $c_2\geq6$ and $L^2>6$.
Now we apply 2.1.3(b). The same arguments as in the previous step show that
the number of irreducible components of $\overline M(L,c_2)$ is at most
that of $\overline M_{\kl_0(-2D)}(\kl_0(-2D),c_2-1)$.
Since
$c_2-1-\frac{(\kl_0(-2D))^2}{2}-3=e+1$, we can use induction over $-e$ and step
1 to show the irreducibility in this case.\\
{\it 4th step:} It remains to consider the case $L^2=2$. Here we apply 2.1.2
with
$d=5$. As above we conclude that the number of irreducible components
of $\overline M(L,c_2)$ is at most that of $\overline
M_{\kl_0(2D)}(\kl_0(2D),c_2+3)$. Since $(\kl_0(2D))^2=L^2+20-8=14$ we  can
conclude by
step 2 or 3.
\qed

Mukai seems to know that all moduli spaces of rank two bundles on
a K3 surface are irreducible (\cite{Mu2}, p.\ 157). Since we could not find a
proof
of this  in
the literature we decided to include the above lemma.

Let $X$ be a K3 surface and $L$ a line bundle on $X$. For any $c_2$ there
exists a coarse moduli space $\overline{M}_s(L,c_2)$ of simple
sheaves of rank two with determinant $L$ and second Chern class $c_2$.
$\overline{M}_s(L,c_2)$ can be realized as a non-separated algebraic space
(\cite{AK}, \cite{KO}).
For any polarization $H$ such that $H-$semi-stabilty implies $H-$stability the
projective manifold $\overline{M}_H(L,c_2)$ is an open subset of
$\overline{M}_s(L,c_2)$.

Note that in the case that $\Pic(X)=\IZ\cdot L$ and $H=L$ any simple vector
bundle
is in fact slope stable. For sheaves the situation is more complicated.

Now let $(\kx,\kl)\lra S$ be a family of K3 surfaces with a line bundle $\kl$
on $\kx$ over a smooth curve $S$. By \cite{AK}, \cite{KO} there exists a
relative moduli space of
simple sheaves, i.e. there exists an algebraic space
$\overline{\km}_s(\kl,c_2)$ and a morphism from it to $S$ such that
the fibre over a point $t\in S$ is isomorphic to $\overline{M}_s(\kl_t,c_2)$.
By a result of Mukai the fibres are smooth \cite{Mu1}. Lemma \ref{irr}
shows that for a family $(\kx,\kl)\lra S$ both
$\overline\km_s(\kl,c_2)$ and $\overline\km_s(\kl,c_2)\lra S$
are smooth (at least over the locally free sheaves).\\
For the following we want to assume that $\Pic(\kx_t)\cong\IZ\cdot\kl_t$ for
$t\not=0$ and $\kl_t^2>0$.
To shorten notation we denote by $Z^*\lra S^*$ the restriction of a family
$Z\lra S$ to $S^*:=S\setminus\{0\}$.

\begin{proposition}\label{defofmod}
Assume that $\overline M_{\kl_t}(\kl_t,c_2)$ is irreducible for $t\not=0$.
Then for any generic ample $H\in \Pic(\kx_0)$
there exists a smooth proper family $Z\lra S$ of
projective manifolds such that $Z^*\lra S^*$ has fibres
$\overline{M}_{\kl_t}(\kl_t,c_2)$ and the fibre over $0$ is isomorphic
to $\overline{M}_{H}(\kl_0,c_2)$. (``The moduli spaces for different $H$ cannot
be separated'')
\end{proposition}
{\it Proof:}
By $\overline{\km}(\kl,c_2)^*\to S^*$ we denote the family of the moduli
spaces $\overline{M}(\kl_t,c_2)$. It is proper over $S^*$ and the fibres
are smooth and
irreducible.\\
{\it Claim:} If $[E]\in \overline M_s(\kl_0,c_2)$ is a point in the
closure $T_0$
of $\overline\km_s(\kl,c_2)^*\setminus\overline{\km}(\kl,c_2)^*$ in
$\overline\km_s(\kl,c_2)$, then $E$ is not semi-stable with respect
to any polarization
$H$: Semi-continuity shows that a point $E$ in the closure has
a subsheaf of rank one with determinant $ \kl_0^{\otimes n}$ with $n>0$. Hence
it is
not semi-stable with respect to any polarization.\\
The set $T_1$ of simple sheaves $[E]\in\overline M_s(\kl_0,c_2)$ which are
not stable with respect to $H$ is a closed subset of
$\overline\km_s(\kl,c_2)$.
We define $Z$ to be the complement of the union of $T_0$ and $T_1$ in
$\overline\km_s(\kl,c_2)$. It is an open subset of $\overline{\km}_s(\kl,c_2)$.
The fibres meet the requirements of the assertion.\\
{\it Claim:} $Z$ is separated:
Any simple sheaf on any of the fibres $\kx_t$ can also be regarded as
a simple coherent sheaf on the complex space $\kx$. Thus $Z$ is a subspace of
the space of all simple sheaves on $\kx$. In order to show that two
points are separated in $Z$ it suffices to separate them in the  bigger space.
Now we apply the criterion of \cite{KO} which says that if two simple coherent
sheaves are not separated then there exists a non-trivial homomorphism between
them. Since any two sheaves parametrized by $Z$ are either supported
on different fibres or stable with respect to the same polarization,
this is excluded.\\
Thus $Z$ is a separated with compact irreducible
fibres over $S^*$. Take
a locally free $E\in \overline M_H(L,c_2)$ and consider a neighbourhood of it
in $\overline{\km}_s(\kl,c_2)$. By the arguments above this
neighbourhood contains locally free simple sheaves
on all the nearby fibres. Hence we can assume that all these sheaves on
$\kx_{t\not=0}$ are stable, since
$\Pic(\kx_t)=\IZ\cdot\kl_t$ for $t\not=0$. This
implies the connectedness of $Z$.
Thus $Z\lra S$ is proper and smooth.\qed

{\bf Proof of the theorem:}
{\it i)} We first show that the result of section 1 generalizes to the case
where
we drop the assumption that $L$ generates $\Pic(X)$. This is done as follows.
By applying \ref{defofmod} to a
deformation of the type 2.1.1 one sees that
$\overline M_{\kl_t}(\kl_t,\frac{\kl_t^2}{2}+3)$ is a deformation of
$\overline M_H(L,\frac{L^2}{2}+3)$ for generic $H$. Since Hodge numbers
are invariant under deformations, both spaces have the same Hodge numbers.
Those of the second were compared in section 1
with the Hodge numbers of the appropriate Hilbert scheme.
By the same trick we can always reduce to the case where the Picard group
is generated by $L$, in particular we can assume that $L$ is ample. \\
{\it ii)}
By applying \ref{defofmod}
to a deformation of type 2.1.2, 2.1.3(a) or 2.1.3(b) we see that $\overline
M(L,c_2)$
is a deformation of $\overline M_H(\kl_0,c_2)$ for generic $H$.
Since $\mu-$stabilty does not change under twisting by line bundles
we have
$\overline M_H(\kl_0,c_2)\cong \overline M_H(\kl_0(2D),c_2+\kl_0.D+D^2)$
(or $\overline M_H(\kl_0(-2D),c_2-\kl_0.D+D^2)$ in case 2.1.3(b)).
The proof of lemma \ref{irr} shows that by applying 2.1.2, 2.1.3(a) and
2.1.3(b) repeatedly we can reduce to the situation of {\it i)}, i.e.
$c_2=\frac{L^2}{2}+3$.
\qed

\begin{corollary}
Let $X$ be an arbitrary K3 surface and $L$ a primitive big and nef line
bundle. As long as a polarization $H$ does not lie on any wall, all
deformation invariants, e.g. Hodge- and Betti numbers,
of $\overline{M}_H(L,c_2)$ are independent of $H$.\qed
\end{corollary}
For similar results compare \cite{G}.
{\footnotesize 
\footnotesize{Max-Planck-Institut f\"ur Mathematik\\
Gottfried-Claren-Str. 26\\
53225 Bonn\\
Germany\\
lothar@mpim-bonn.mpg.de\\
huybrech@mpim-bonn.mpg.de}
\end{document}